\documentclass{amsart}
\usepackage{amsmath}
\usepackage{amssymb}
\usepackage{amscd}
\usepackage[all]{xy}
\usepackage{color}
\usepackage{fancyhdr}
\usepackage{mathrsfs}

\theoremstyle{plain}
\newtheorem{theorem}{Theorem}[section]
\newtheorem{proposition}[theorem]{Proposition}
\newtheorem{lemma}[theorem]{Lemma}
\newtheorem{corollary}[theorem]{Corollary}
\newtheorem{definition}[theorem]{Definition}

\theoremstyle{question}

\newtheorem{example}[theorem]{Example}

\theoremstyle{remark}
\newtheorem{remark}[theorem]{Remark}

\def\bc{\begin{center}}
\def\ec{\end{center}}

\def\Spec{{\rm Spec}}

\def\ann{{\rm ann}}

\def\Ass{{\rm Ass}}

\def\height{{\rm ht}}

\def\GV{{\rm GV}}
\def\Spec{{\rm Spec}}

\def\Ass{{\rm Ass}}
\def\Spec{{\rm Spec}}

\def\wmax{w\mbox{-}\max}
\def\tmax{t\mbox{-}\max}
\def\starmax{\ast\mbox{-}\max}

\begin{document}
\title[Two generalizations of Krull domains] {Two generalizations of Krull domains}

\author [Xing] {Shiqi Xing$^{\ast}$}
\address{(Xing) College of Applied Mathematics, Chengdu University of Information Technology, Chengdu, Sichuan 610225, China}
\email{sqxing@yeah.net}
\author [Anderson] {D. D. Anderson}
\address{(Anderson) Department of Mathematics, The University of Iowa, Iowa City, IA 52242-1419, USA}
\email{dan-anderson@uiowa.edu}
\author [Zafrullah] {Muhammad Zafrullah}
\address{(Zafrullah) Department of Mathematics, Idaho State University, Pocatello, ID 83209-8085, USA}
\email{mzafrullah@usa.net}


\thanks{$2010$ Mathematics Subject Classification: 13A15, 13A18, 13F05, 13G05}

\thanks{$\ast$ Corresponding Author}

\date{\today}

\begin{abstract}
In this paper we introduce two new generalizations of Krull domains: $\ast$-almost independent rings of Krull type ($\ast$-almost IRKTs) and $\ast$-almost generalized Krull domains ($\ast$-AGKDs), neither of which need be integrally closed. We characterize them using certain types of $\ast$-homogeneous ideals. To do this we introduce $\ast$-almost
super-homogeneous ideals and $\ast$-almost super-SH domains.  We prove that a domain $D$ is a $\ast$-almost IRKT if and only if $D$ is a $\ast$-almost super-SH domain and that a domain is a $\ast$-AGKD if and only if $D$ is a type 1 $\ast$-almost super-SH domain. Further, we study $\ast$-almost factorial general-SH domains ($\ast$-afg SH domains) and we prove that a domain $D$ is a $\ast$-afg-SH domain if and only if $D$ is a $\ast$-IRKT and an AGCD-domain.
\end{abstract}

\maketitle

\section{Introduction}
It is well-known that Krull domains play a central role in the development of multiplicative ideal theory. The concept of a Krull domain has been generalized in many different ways, for example, by independent rings of Krull type, generalized Krull domains and weakly Krull domains. There is an important commonness in the above domains, i.e., they are all $\mathscr{F}$-IFC domains. Recall that a set $\mathscr{F}$ of prime ideals in a domain $D$ is a \emph{defining family} if $D=\bigcap\{D_{P}\mid P\in \mathscr{F}\}$. Further, $\mathscr{F}$ is \emph{of finite character} (or \emph{locally finite}) if every nonzero nonunit of $D$ belongs to at most finitely many members of $\mathscr{F}$ and $\mathscr{F}$ is \emph{independent} if no two members of $\mathscr{F}$ contain a common nonzero prime ideal. As in \cite{AZ03}, a domain $D$ is called an \emph{$\mathscr{F}$-IFC domain} if $D$ has a defining family $\mathscr{F}$ such that $\mathscr{F}$ is independent and of finite character. Now suppose that $D$ is a $\mathscr{F}$-IFC domain. Then $D$ is called a \emph{weakly Krull domain} (WKD) in \cite{AHZ} if $X^{(1)}(D)=\mathscr{F}$, where $X^{(1)}(D)$ is the set of height-one prime ideals of $D$. We can further put conditions on $D_{P}$ for $P\in \mathscr{F}$. If $D$ is an $\mathscr{F}$-IFC domain and $D_{P}$ is a valuation domain for each $P\in\mathscr{F}$, then we get the \emph{independent rings of Krull type} (IRKTs) of Griffin \cite{Gr}. If $D$ is a WKD and $D_{P}$ is a valuation domain for each $P\in\mathscr{F}$, then we get the \emph{generalized Krull domains} (GKDs) of Ribenboim \cite{R}. In particular, if $D$ is a WKD and $D_{P}$ is a DVR for each $P\in\mathscr{F}$, then $D$ is precisely a \emph{Krull domain}. Our original motivation for this paper was to give two classes of generalizations of Krull domains by $\mathscr{F}$-IFC domains.

Let us denote the set of positive integers by $\mathbb{N}$. Recall from \cite{AZ02} that a domain $D$ is called an \emph{almost valuation domain} (AV-domain) if for $0\neq a,b\in D$, there is an $n=n(a,b)\in \mathbb{N}$ with $a^{n}\mid b^{n}$ or $b^{n}\mid a^{n}$. It is clear that every valuation domain is an AV-domain. Now AV-domains are of interest in that by using AV-domains, many classical results on valuation domains can be extended to the general theory of almost factoriality. For example, a domain $D$ is called an \emph{almost Pr\"{u}fer domain} (\emph{AP-domain}) in \cite{AZ02} if for any $0\neq a,b\in D$, there is an  $n=n(a,b)\in \mathbb{N}$ with $(a^{n}, b^{n})$ invertible. It is shown in \cite[Theorem 5.8]{AZ02} that a domain $D$ is an AP domain if and only if $D_{P}$ is an AV-domain for each maximal ideal $P$ of $D$. Also, a domain $D$ is called an \emph{almost Pr\"{u}fer $v$-multiplication domain} (\emph{AP$v$MD}) in \cite{L02} if for $0\neq a,b\in D$, there is an $n=n(a, b)\in \mathbb{N}$ with $(a^{n}, b^{n})$ $t$-invertible. It is shown in \cite[Theorem 2.3]{L02} that a domain $D$ is an AP$v$MD if and only if $D_{P}$ is an AV-domain for each maximal $t$-ideal $P$ of $D$. Now in the definition of an IRKT, we can use an AV-domain instead of a valuation domain to define an almost independent ring of Krull type. If $D$ is a $\mathscr{F}$-IFC domain and $D_{P}$ is an AV-domain for any $P\in \mathscr{F}$, then $D$ is said to be an \emph{almost independent ring of Krull type} (\emph{almost IRKT}). Accordingly, if $D$ is a WKD and $D_{P}$ is an AV-domain for any $P\in \mathscr{F}$, $D$ is said to be an \emph{almost generalized Krull domain} (\emph{AGKD}). In this paper, we shall create a suitable theory of unique factorization of ideals as in \cite{AZ} and study them in a slightly more general setting using finite character star-operations.

Let $\ast$ be a finite character star-operation on a domain $D$. Denote the set of maximal $\ast$-ideals by $\starmax(D)$. Then $D=\bigcap\{D_{P}\mid P\in \starmax(D)\}$ by \cite[Theorem 7.2.11]{WK02}. Hence $\starmax(D)$ is a defining family on $D$. As in \cite{AZ}, a $\starmax(D)$-IFC domain is said to be \emph{$\ast$-h-local}. Indeed, following \cite[page 136]{FS}, a $\ast$-h-local domain can be also called $h_{\mathcal{P}}$-local where $\mathcal{P}=\starmax(D)$. In particular, a $d$-h-local domain is precisely a \emph{h-local domain} of Matlis \cite{M} and a $t$-h-local domain is precisely a \emph{$h_{\mathcal{U}}$-local domain} \cite {BKW}. Using $\ast$-h-local domains,  WKDs, IRKTs, GKDs, and Krull domains are redefined in \cite{AZ} as $\ast$-WKDs, $\ast$-IRKTs, $\ast$-GKDs, and $\ast$-Krull domain, respectively. More precisely, a domain $D$ is called a \emph{$\ast$-WKD} if $D$ is $\ast$-h-local and $\starmax(D)=X^{(1)}(D)$; a domain $D$ is called a \emph{$\ast$-IRKT} if $D$ is $\ast$-h-local and $D_{P}$ is a valuation domain for each $P\in \starmax(D)$; a domain $D$ is called a \emph{$\ast$-GKD} if $D$ is a $\ast$-WKD and $D_{P}$ is a valuation domain for each $P\in X^{(1)}(D)$; a domain is called
a \emph{$\ast$-Krull domain} if $D$ is a $\ast$-WKD and $D_{P}$ is a DVR for each $P\in \starmax(D)$. It is easy to check that a $t$-WKD (resp., $t$-IRKT, $t$-GKD and $t$-Krull domain) is just a WKD (resp., IRKT, GKD and Krull domain) while a $d$-WKD (resp., $d$-IRKT, $d$-GKD and $d$-Krull domain) is a one-dimensional finite character domain (resp., finite character Pr\"{u}fer domain, one dimensional finite character Pr\"{u}fer domain and a Dedekind domain). Using $\ast$-homogeneous ideals, the second and third authors have given some nice characterizations for these domains. Let $\ast$ be a finite character star-operation on a domain $D$. Recall from \cite{AZ} that a nonzero ideal $A$ of $D$ is called \emph{$\ast$-homogeneous} if $A$ is finitely generated and $A$ is contained in a unique maximal $\ast$-ideal. The unique maximal $\ast$-ideal containing $A$ is often denoted by $M(A)$. It is shown in \cite[Theorem 4]{AZ} that a domain $D$ is a $\ast$-h-local domain if and only if $D$ is a $\ast$-SH domain, where $D$ is called a \emph{$\ast$-semi-homogeneous domain} (\emph{$\ast$-SH domain}) if every proper nonzero principal ideal of $D$ is a $\ast$-product of $\ast$-homogeneous ideals. Furthermore, if $A$ is $\ast$-homogeneous and $M(A)=\sqrt{A_{\ast}}$, then $A$ is called a \emph{type 1 $\ast$-homogeneous} ideal; if $A$ is $\ast$-homogenous and $A_{\ast}=(M(A)^{n})_{\ast}$ for some $n\geq 1$, then $A$ is called a \emph{type 2 $\ast$-homogeneous ideal}. Accordingly,  a domain $D$ is called a \emph{type 1 $\ast$-SH domain }if every proper nonzero principal ideal of $D$ is a finite $\ast$-product of type 1 $\ast$-homogeneous ideals of $D$; a domain is called a \emph{type 2 $\ast$-SH domain} if every proper nonzero principal ideal of $D$ is a finite $\ast$-product of type 2 $\ast$-homogeneous ideals of $D$. It is shown in \cite[Theorem 7]{AZ} that a domain $D$ is a $\ast$-WKD if and only if $D$ is a type 1 $\ast$-SH domain. It is shown in \cite[Theorem 8]{AZ} that a domain $D$ is a $\ast$-Krull domain if and only if $D$ is a type 2 $\ast$-SH domain. Also, a $\ast$-homogeneous ideal $A$ of $D$ is called \emph{$\ast$-super homogeneous} if each $\ast$-homogeneous ideal containing $A$ is $\ast$-invertible and a domain $D$ is called a \emph{$\ast$-super-SH domain} if every proper nonzero principal ideal of $D$ is a finite $\ast$-product of $\ast$-super-homogeneous ideals of $D$. It is shown
\cite[Theorem 10]{AZ} that a domain $D$ is a $\ast$-IRKT if and only if $D$ is a $\ast$-super-SH domain. In analogy with $\ast$-super-homogeneous ideals, we introduce $\ast$-almost super-homogeneous ideals and $\ast$-almost super-SH domains in Section 2. Here a
\emph{$\ast$-almost super-homogeneous ideal} $A$ is a $\ast$-invertible $P$-$\ast$-homogeneous ideal with the additional condition that given $b_{1}, \dots,b_{s}\in P$ with $A^{r}\subseteq (b_{1},\dots, b_{s})_{\ast}$ for some $r\in \mathbb{N}$, there exists an $n\in \mathbb{N}$ with $(b_{1}^{n},\dots,b_{s}^{n})$ $\ast$-invertible. Accordingly, a domain $D$ is called a \emph{$\ast$-almost super-SH domain} if every nonzero proper principal ideal of $D$ is a $\ast$-product of $\ast$-almost super-homogeneous ideals. This fills a gap left in \cite{AZ}. The $\ast$-almost factorial-SH domains introduced in \cite{AZ} are integrally closed while the $\ast$-almost super-SH domains introduced in this paper need not be integrally closed. In Section 3, we study $\ast$-almost IRKTs and we prove in Theorem \ref{011} that a domain $D$ is a $\ast$-almost IRKT if and only if $D$ is a $\ast$-almost super SH domain. In Section 4, we study $\ast$-AGKDs and we prove in Theorem \ref{014} that a domain $D$ is a $\ast$-AGKD if and only if $D$ is a type 1 $\ast$-almost super-SH domain. Furthermore, in Section 5, we study $\ast$-almost factorial general-SH domains ($\ast$-afg-SH domains) and we prove in Theorem \ref{019} that a domain $D$ is a $\ast$-afg-SH domain if and only if $D$ is an AGCD-domain and a $\ast$-almost IRKT, where a domain $D$ is called an \emph{almost GCD domain} (\emph{AGCD-domain}) in \cite{Z} if for $0\neq a,b\in D$, there exists an $n=n(a,b)\in \mathbb{N}$ with $(a^{n}, b^{n})_{v}$ principal.

As our work involves star-operations, we provide a quick review.
Let $D$ be a domain with quotient field $K$ and let $F(D)$ be the set of nonzero fractional ideals
of $D$. A \textit{star-operation} on $D$ is a map $\ast : F(D)\to F(D)$ such that for all $A, B\in F(D)$ and $0\neq x\in K$

(1)\ \ $(x)_{\ast}=(x)$ and $(xA)_{\ast}=xA_{\ast}$,

(2)\ \ $A\subseteq A_{\ast}$, $A_{\ast}\subseteq B_{\ast}$ whenever $A\subseteq B$, and

(3)\ \ $(A_{\ast})_{\ast}=A_{\ast}$.

We note that for $A,B\in F(D)$, $(AB)_{\ast }=(A_{\ast}B)_{\ast}=(A_{\ast }B_{\ast })_{\ast }$ and call it the \emph{$\ast $-product}. A fractional ideal $A$ is called a $\ast $\emph{-fractional ideal} if $A=A_{\ast }$ and $A$ is called a fractional ideal of \emph{$\ast$-finite type} if there exists a finitely generated fractional ideal $B\in F(D)$ such that $A_{\ast}=B_{\ast }$. A star-operation $\ast $ is said to be \emph{of finite character} or
\emph{ of finite type} if $ A_{\ast }=\bigcup \{B_{\ast }\mid 0\neq B$ is a
finitely generated fractional ideal contained in $A\}$ for each
$A\in F(D)$. For $A\in F(D)$, define $A^{-1}:=\{x\in K\mid
xA\subseteq D\}$ and call $A$ $\ast $\emph{-invertible} if
$(AA^{-1})_{\ast }=D$. If $\ast$ is a star-operation on a domain $D$, then $\ast$ always induces two finite character star-operations, $\ast_{s}$ and $\ast_{w}$. Let $A\in F(D)$. Then $A_{\ast_{s}}=\bigcup\{B_{\ast}\mid 0\neq B$ f.g. and $B\subseteq A\}$, and $A_{\ast_{w}}=\{x\in K\mid xJ\subseteq A $ for some nonzero f.g. ideal $J$ with $J_{\ast}=D\}$. The classical star-operations are the $v$-, $t$-,
$w$-operations. Let $A$ be a nonzero fractional ideal of $D$. Then
$A_v := (A^{-1})^{-1}$, $A_t :=\bigcup \{B_v \mid 0\neq B$ f.g. and
$B\subseteq A\}=A_{v_{s}}$, and $A_w :=\{x
\in K \mid Jx \subseteq A$ for some $J
\in \GV(D) \}=A_{v_{w}}$, where $\GV(D)=\{J\mid J$ is a nonzero f.g. ideal of $D$ with $J^{-1}=D\}$.  We now proceed to state and prove our main
results.

\section{$\ast$-almost super-homogenous ideals}
In this section we introduce $\ast$-almost super-homogeneous ideals and $\ast$-almost super-SH domains.
Suppose that $A$ is a $\ast$-homogeneous ideal of a domain $D$. If $P$ is the unique maximal $\ast$-ideal containing $A$, then $A$ is said to be \emph{$P$-$\ast$-homogeneous}. If both $A$ and $B$ are $P$-$\ast$-homogeneous, we say that $A$ is \emph{similar} to $B$, denoted by $A\sim B$. Now we start by the following definition.
\begin{definition}\label{001}
\rm{Let $\ast$ be a finite character star-operation on a domain $D$.
\begin{itemize}
\item[(1)]An ideal $A$ of $D$ is called \emph{$\ast$-almost super-homogeneous} if
\begin{itemize}
\item[(i)] $A$ is a $\ast$-invertible $P$-$\ast$-homogeneous ideal, and
\item[(ii)] given $b_{1}, \dots,b_{s}\in P$ with $A^{r}\subseteq (b_{1},\dots, b_{s})_{\ast}$ for some $r\in \mathbb{N}$, there exists an $n\in \mathbb{N}$ with $(b_{1}^{n},\dots,b_{s}^{n})$ $\ast$-invertible.
\end{itemize}
\item[(2)] $D$ is called a \emph{$\ast$-almost super-SH domain} if every nonzero proper principal ideal of $D$ is a $\ast$-product of $\ast$-almost super homogeneous ideals.
\end{itemize}}
\end{definition}

\remark\label{23}
In Definition \ref {001} (1), the $n$ depends on $b_{1},\dots, b_{s}$.\\

Next we investigate the properties of $\ast$-almost super-homogeneous ideals. We need the following lemmas.

\begin{lemma}\label{002}
Let $A=(a_{1},\dots, a_{k})$ be an ideal of a domain $D$. Then $A^{nk}\subseteq (a_{1}^{n},\dots,a_{k}^{n})\subseteq A^{n}$ for any $n\in \mathbb{N}$.
\end{lemma}
\begin{proof}
It is clear that $(a_{1}^{n}, \dots, a_{k}^{n})\subseteq A^{n}$. Let $a_{1}^{l_{1}}\cdots a_{k}^{l_{k}}\in A^{nk}$, $\sum_{i=1}^{k}l_{k}=nk$. Then there is some $l_{j}\geq n$. Otherwise, as $l_{i}<n$ for all $i$, we have $\sum_{i=1}^{k}l_{k}<kn$, a contradiction. Hence $a_{1}^{l_{1}}\cdots a_{k}^{l_{k}}\in (a_{1}^{n},\dots,a_{k}^{n})$. It follows that $A^{nk}\subseteq (a_{1}^{n},\dots,a_{k}^{n})$. So $A^{nk}\subseteq (a_{1}^{n}, \dots, a_{k}^{n})\subseteq A^{n}$.
\end{proof}

\begin{lemma}\label{003}
Let $\ast$ be a finite character star-operation on a domain $D$ and $\{a_{\alpha}\}\subseteq D\setminus \{0\}$. If $(\{a_{\alpha}\})_{\ast}$ is $\ast$-invertible, then $(\{a_{\alpha}^{n}\})_{\ast}=((\{a_{\alpha}\})^{n})_{\ast}$ for any $n\in \mathbb{N}$.
\end{lemma}
\begin{proof}
See \cite[Lemma 2.2]{L}.
\end{proof}

\begin{proposition}\label{005}
Let $\ast$ be a finite character star-operation on a domain $D$ and $A$ a $P$-$\ast$-almost super-homogeneous ideal of $D$.
\begin{itemize}

\item[(1)] If $(b_{1},\dots, b_{s})$ is a $P$-$\ast$-homogeneous ideal of $D$, then $(A^{n}+(b_{1}^{n},\dots,b_{s}^{n}))_{\ast}=(A^{n})_{\ast}$ or $(A^{n}+(b_{1}^{n},\dots,b_{s}^{n}))_{\ast}=(b_{1}^{n},\dots,b_{s}^{n})_{\ast}$ for some $n\in \mathbb{N}$.

\item[(2)] If $B$ is a $P$-$\ast$-almost super-homogeneous ideal of $D$, then there exists an $n\in \mathbb{N}$ such that  $B^{n}\subseteq (A^{n})_{\ast}$ or $A^{n}\subseteq (B^{n})_{\ast}$.
\item[(3)] If $B$ is a $P$-$\ast$-almost super-homogeneous ideal of $D$, then so is $AB$.
\item[(4)] $A^{n}$ is $\ast$-almost super-homogeneous for any positive integer $n$.
\end{itemize}
\end{proposition}
\begin{proof}
$(1)$ \ \ Let $A=(a_{1}, \dots, a_{k})$. Then $A\subseteq (a_{1},\dots, a_{k}, b_{1},\dots, b_{s})_{\ast}\subseteq P$. Hence for some $n\in \mathbb{N}$, $(a_{1}^{n},\dots, a_{k}^{n}, b_{1}^{n},\dots, b_{s}^{n})$ is $\ast$-invertible. So $(((a_{1}^{n},\dots, a_{k}^{n})+(b_{1}^{n},\dots, b_{s}^{n}))$ $(a_{1}^{n},\dots, a_{k}^{n}, b_{1}^{n},\dots, b_{s}^{n})^{-1})_{\ast}=((a_{1}^{n}\dots,a_{k}^{n},b_{1}^{n},\dots,b_{s}^{n})(a_{1}^{n},\dots, a_{k}^{n}, b_{1}^{n}, \dots ,  b_{s}^{n})^{-1})_{\ast}=D$. It follows that $(a_{1}^{n},\dots,a_{k}^{n})(a_{1}^{n}\dots,a_{k}^{n},b_{1}^{n},\dots,b_{s}^{n})^{-1}\nsubseteq P$ or $(b_{1}^{n},\dots, b_{s}^{n})(a_{1}^{n}\dots,a_{k}^{n} $ $,b_{1}^{n},\dots,b_{s}^{n})^{-1}\nsubseteq P$. We claim that $(a_{1}^{n},\dots,a_{k}^{n})(a_{1}^{n},\dots,a_{k}^{n},b_{1}^{n},\dots,b_{s}^{n})^{-1}$ and $(b_{1}^{n},\dots $ $, b_{s}^{n})(a_{1}^{n}, \dots,a_{k}^{n},b_{1}^{n},\dots,b_{s}^{n})^{-1}$ can not be contained in any maximal $\ast$-ideal other than $P$. In fact, if $(b_{1}^{n},\dots,b_{s}^{n})(a_{1}^{n},\dots,a_{k}^{n},b_{1}^{n},\dots,b_{s}^{n})^{-1}\subseteq Q$ for some $Q\in\starmax(D)$, then $(b_{1}^{n},\dots,b_{s}^{n})_{\ast}\subseteq (Q(a_{1}^{n},\dots,a_{k}^{n},b_{1}^{n},\dots,b_{s}^{n}))_{\ast}\subseteq Q_{\ast}=Q$. Hence $b_{1},\dots, b_{s}\in Q$. So $Q=P$ since $(b_{1},\dots,b_{s})$ is $P$-homogeneous. Similarly we can show that $(a_{1}^{n},\dots,a_{k}^{n})(a_{1}^{n},\dots,a_{k}^{n},b_{1}^{n},\dots,b_{s}^{n})^{-1}$ can not be contained in any maximal $\ast$-ideal other than $P$.
Thus $((a_{1}^{n},\dots,a_{k}^{n})(a_{1}^{n}\dots,a_{k}^{n},b_{1}^{n},\dots,b_{s}^{n})^{-1})_{\ast}=D$ or $((b_{1}^{n},\dots, b_{s}^{n})(a_{1}^{n},$
$\dots,a_{k}^{n},b_{1}^{n},\dots,b_{s}^{n})^{-1})_{\ast}=D$. In the first case, we have $(A^{n})_{\ast}=(a_{1}^{n},\dots,a_{k}^{n})_{\ast} =(a_{1}^{n}\dots,a_{k}^{n},b_{1}^{n},\dots,b_{s}^{n})_{\ast}$. So $(A^{n})_{\ast}=(A^{n}+(b_{1}^{n},\dots, b_{s}^{n}))_{\ast}$. In the second case, we have
$(b_{1}^{n},\dots, b_{s}^{n})_{\ast}=(a_{1}^{n},\dots,a_{k}^{n},b_{1}^{n},\dots,b_{s}^{n})_{\ast}$. So $(b_{1}^{n},\dots, b_{s}^{n})_{\ast}=(A^{n}+(b_{1}^{n},\dots, b_{s}^{n}))_{\ast}$.

$(2)$\ \ Suppose that $B=(b_{1},\dots, b_{s})$. Now by (1), for some $n\in \mathbb{N}$, either $(b_{1}^{n}\dots,b_{s}^{n})_{\ast}\subseteq (A^{n}+(b_{1}^{n},\dots, b_{s}^{n}))_{\ast}=(A^{n})_{\ast}$ or $(A^{n})_{\ast}\subseteq (A^{n}+(b_{1}^{n},\dots, b_{s}^{n}))_{\ast}=(b_{1}^{n},\dots,b_{s}^{n})_{\ast}$. Since $B$ is $\ast$-invertible, $(B^{n})_{\ast}=(b_{1}^{n},\dots, b_{s}^{n})_{\ast}$ by Lemma \ref{003}. So $B^{n}\subseteq (A^{n})_{\ast}$ or $A^{n}\subseteq (B^{n})_{\ast}$.

$(3)$ \ \ By \cite[Proposition 2]{AZ}, it follows that $AB$ is $P$-$\ast$-homogeneous and certainly $AB$ is $\ast$-invertible. Let $C=(c_{1},\dots, c_{l})$ be a $P$-$\ast$-homogeneous ideal with $(AB)^{r}\subseteq C_{\ast}$ for some $r\in \mathbb{N}$. By (3) we have $A^{n}\subseteq (B^{n})_{\ast}$ or $B^{n}\subseteq (A^{n})_{\ast}$ for some $n\in \mathbb{N}$. Without loss of generality, suppose that $A^{n}\subseteq (B^{n})_{\ast}$. Then by Lemma \ref{002}, $A^{2nlr}\subseteq (A^{2nlr})_{\ast}\subseteq ((AB)^{nlr})_{\ast}\subseteq ((c_{1},\dots, c_{l})^{nl})_{\ast}\subseteq (c_{1}^{n},\dots, c_{l}^{n})_{\ast}\subseteq ((c_{1},\dots,c_{l})^{n})_{\ast}\subseteq P$. Since $A$ is $P$-$\ast$-almost super-homogeneous, there exists some $m\in\mathbb{N}$ with $(c_{1}^{mn},\dots, c_{l}^{mn})$ $\ast$-invertible. Hence $AB$ is $P$-$\ast$-almost super-homogeneous.

$(4)$\ \ This follows from $(3)$.
\end{proof}

\begin{corollary}\label{006}
Let $\ast$ be a finite character star-operation on a domain $D$ and $A$ a $P$-$\ast$-almost super-homogeneous ideal of $D$. If $B$ is a $\ast$-invertible and $A^{r}\subseteq B_{\ast}\neq D$ for some $r\in \mathbb{N}$, then $B$ is $P$-$\ast$-almost super-homogeneous.
\end{corollary}
\begin{proof} By Proposition \ref{005} (4), it follows that $A^{r}$ is $P$-$\ast$-almost super-homogeneous. Let $B^{l}\subseteq (c_{1},\dots,c_{k})_{\ast}\subseteq P$ for some $l\in\mathbb{N}$. Then $A^{rl}\subseteq (B^{l})_{\ast}\subseteq (c_{1},\dots,c_{k})_{\ast}$. Hence, there exists some $n\in \mathbb{N}$ with $(c_{1}^{n},\dots,c_{k}^{n})$ $\ast$-invertible. So $B$ is also $P$-$\ast$-almost super-homogeneous.
\end{proof}
\begin{proposition}\label{007}
Let $\ast$ be a finite character star-operation on a domain $D$. If $A$ is a $\ast$-super-homogeneous ideal of $D$, then $A$ is a $\ast$-almost
super-homogeneous ideal.
\end{proposition}
\begin{proof} Suppose that $A^{r}\subseteq (c_{1},\dots,c_{k})_{\ast}\subseteq P$ for some $r\in\mathbb{N}$. Now as $A^{r}$ is $\ast$-super-homogeneous and similar to $A$ by \cite[Theorem 1.11]{HZ}, $(c_{1},\dots,c_{k})$ is $\ast$-invertible. Hence $A$ is $M(A)$-$\ast$-almost super-homogeneous.
\end{proof}

\begin{corollary}\label{008}
If $D$ is a $\ast$-super-SH domain, then $D$ is a $\ast$-almost super-SH domain.
\end{corollary}
\begin{proof}
This follows from Proposition \ref{007}.
\end{proof}
We now give a uniqueness result for $\ast$-products of $\ast$-almost homogeneous ideals.
\begin{theorem}\label{009}
Let $\ast$ be a finite character star-operation on a domain $D$ and let $A_{1},\dots, A_{n}$ be $\ast$-almost super-homogeneous ideals of $D$. Then the $\ast$-product $(A_{1}\cdots A_{n})_{\ast}$ can be expressed uniquely, up to order, as a product of pairwise $\ast$-comaximal $\ast$-almost super-homogeneous ideals.
\end{theorem}
\begin{proof} Write $A=(A_{1}\cdots A_{n})_{\ast}$. Let $M(A_{i_{1}}),\dots, M(A_{i_{s}})$ be the distinct maximal $\ast$-ideals among $M(A_{1}), \dots, M(A_{n})$. Set $B_{k}:=\prod\{A_{j}\mid A_{j}\sim A_{i_{k}}\}$ ($k=1,\dots, s$). Then by Proposition \ref{005} (3), $B_{1},\dots, B_{s}$ are $\ast$-almost
super-homogeneous ideals of $D$ that are pairwise $\ast$-comaximal and $A=(B_{1}\cdots B_{s})_{\ast}$. The uniqueness follows from \cite[Theorem 3]{AZ}.

\end{proof}

\section{$\ast$-almost IRKTs}
In this section we introduce $\ast$-almost IRKTs. Using AV-domains instead of valuation domains, we define a $\ast$-almost IKRT as follows.

\begin{definition}\label{010}
\rm{A domain $D$ is called a \emph{$\ast$-almost independent ring of Krull type} (\emph{$\ast$-almost IRKT}) if $D$ is a $\ast$-h-local domain and $D_{P}$ is an AV-domain for each $P\in\starmax(D)$.}
\end{definition}

Recall from \cite{L} that a domain $D$ is called an \emph{almost Pr\"{u}fer $\ast$ multiplicatiion domain} (\emph{AP$\ast$MD}) if for $0\neq a, b \in D$, there exists an $n=n(a,b)\in \mathbb{N}$ with $(a^{n},b^{n})$ $\ast_{s}$-invertible. Let $\ast$ a finite character star-operation on a domain $D$. Then $D$ is an AP$\ast$MD if and only if $D_{P}$ is an AV-domain for each $P\in \starmax(D)$ \cite[Theorem 2.4]{L}. Hence it is clear that $D$ is a $\ast$-almost IRKT if and only if $D$ is a $\ast$-h-local domain and an AP$\ast$MD. Next we prove that $D$ is a $\ast$-almost IRKT if and only if $D$ is a $\ast$-almost super-SH domain.

\begin{theorem}\label{011}
 Let $D$ be a domain and $\ast$ a finite character star-operation on $D$. Then the following statements are equivalent for $D$.
\begin{itemize}
\item[(1)] $D$ is a $\ast$-almost IRKT.
\item[(2)] $D$ is a $\ast$-h-local domain and an AP$\ast$MD.
\item[(3)] $D$ is a $\ast$-h-local domain and every $\ast$-invertible $\ast$-homogeneous ideal of $D$ is $\ast$-almost super-homogeneous.
\item[(4)] $D$ is a $\ast$-almost super-SH domain.
\end{itemize}
\end{theorem}
\begin{proof}$(1)\Leftrightarrow (2)$ \ \ This follows from \cite[Theorem 2.4]{L}.

$(1)\Rightarrow (3)$ \ \ Suppose that $A$ is $\ast$-invertible and $P$-$\ast$-homogeneous. Let $B=(b_{1},\dots,b_{k})$ be $P$-$\ast$-homogeneous with $A^{r}\subseteq B_{\ast}$. Since $D_{P}$ is an AV-domain, $D_{P}$ is an AB-domain by \cite[Theorem 5.6]{AZ02}. Hence by \cite[Lemma 4.3]{AZ02} there exists some $n\in \mathbb{N}$ such that $(b_{1}^{n},\dots, b_{k}^{n})D_{P}=(b_{1}^{n}/1,\dots,b_{k}^{n}/1)$ is a principal ideal of $D_{P}$. Let $Q\in \starmax(D)$ with $Q\neq P$. Then $B\nsubseteq Q$ since $B$ is $P$-$\ast$-homogeneous. Take $b_{j}\in B\setminus Q$. Then, $b_{j}^{n}\notin Q$. Hence $(b_{1}^{n},\dots, b_{k}^{n})\nsubseteq Q$ and so $(b_{1}^{n},\dots, b_{k}^{n})D_{Q}=D_{Q}$. It follows that $(b_{1}^{n},\dots, b_{k}^{n})D_{M}$ is a locally principal ideal of $D$ for each $M\in \starmax(D)$. So $(b_{1}^{n},\dots, b_{k}^{n})$ is $\ast$-invertible by \cite[Theorem 7.2.15]{WK02}. Consequently, $A$ is $\ast$-almost
super-homogeneous.

$(3)\Rightarrow(4)$\ \ Since $R$ is a $\ast$-h-local domain, it follows from \cite[Theorem 5]{AZ} that $R$ is a $\ast$-SH domain. Let $xD$ be a proper principal ideal of $D$. Then $xD=(A_{1}\cdots A_{k})_{\ast}$ where each $A_{i}$ is $\ast$-homogeneous. Since $xD$ is $\ast$-invertible, each $A_{i}$ is $\ast$-invertible. Hence each $A_{i}$ is $\ast$-almost super-homogeneous. It follows that $D$ is a $\ast$-almost super-SH domain.

$(4)\Rightarrow (1)$\ \ Suppose that $D$ is a $\ast$-almost super-SH domain. Then $D$ is a $\ast$-SH domain. Hence $D$ is a $\ast$-h-local domain by \cite[Theorem 5]{AZ}. We only need to prove that $D_{P}$ is an AV-domain for each $P\in \starmax(D)$. For given $P\in \starmax(D)$, take $0\neq x\in P$. Then by Theorem \ref{009} $xD=(A_{1}\cdots A_{k})_{\ast}$ where the $A_{i}$ are mutually $\ast$-comaximal $\ast$-almost super-homogeneous ideals of $D$. Hence there exists some $A_{i}$ such that $A_{i}\subseteq P$. If $j\neq i$ and $A_{j}\subseteq P$, then $(A_{i}+A_{j})_{\ast}\subseteq P$. But $A_{i}$ and $A_{j}$ are $\ast$-comaximal, so $D=(A_{i}+A_{j})_{\ast}\subseteq P$, which is a contradiction. So there is only one $A_{i}$ such that $A_{i}\subseteq P=M(A_{i})$. Since $A_{i}$ is $\ast$-almost super-homogeneous, $A_{i}$ is $\ast$-invertible. Hence $(A_{i})_{\ast}D_{P}=A_{i}D_{P}$  by \cite[Corollary 7.2.16]{WK02}. It follows that
$(A_{i})_{\ast}=((A_{i})_{\ast})_{\ast_{w}}=\bigcap\{(A_{i})_{\ast}D_{Q}\mid Q\in\starmax(D)\}=A_{i}D_{P}\bigcap D$. So $xD_{P}\bigcap D=(A_{1}\cdots A_{k})_{\ast}D_{P}\bigcap D=(A_{1}\cdots A_{k})D_{P}\bigcap D=A_{i}D_{P}\bigcap D=(A_{i})_{\ast}$. For convenience we write $xD_{P}\bigcap D=A_{\ast}$, where $A$ is $P$-$\ast$-almost super-homogeneous. Let $0\neq y\in P$. Then similarly we get that $yD_{P}\bigcap D=B_{\ast}$, where $B$ is also $P$-$\ast$-almost
super-homogeneous. Thus $A^{n}\subseteq (B^{n})_{\ast}$ or $B^{n}\subseteq(A^{n})_{\ast}$ for some $n\in N$ by Proposition \ref{005}(3). Now we claim that $x^{n}D_{P}\subseteq y^{n}D_{P}$ or $y^{n}D_{P}\subseteq x^{n}D_{P}$. In fact, if $A^{n}\subseteq (B^{n})_{\ast}$, then $((xD_{P}\bigcap D)^{n})_{\ast}=(A^{n})_{\ast}\subseteq (B^{n})_{\ast}=
((yD_{P}\bigcap D)^{n})_{\ast}$. Since $B$ is $\ast$-invertible, $(B^{n})_{\ast}=
((yD_{P}\bigcap D)^{n})_{\ast}$ is $\ast$-invertible. Hence $(yD_{P}\bigcap D)^{n}D_{P
}=((yD_{P}\bigcap D)^{n})_{\ast}D_{P}$ by \cite[Corollary 7.2.16]{WK02}. So we have $x^{n}D_{P}=(xD_{P})^{n}=((xD_{P}\bigcap D)D_{P})^{n}=(xD_{P}\bigcap D)^{n}D_{P}=((xD_{P}\bigcap D)^{n})_{\ast}D_{P}\subseteq ((yD_{P}\bigcap D)^{n})_{\ast}D_{P}=(yD_{P}\bigcap D)^{n}D_{P}=y^{n}D_{P}$. Similarly we can prove that if $B^{n}\subseteq (A^{n})_{\ast}$, then $y^{n}D_{P}\subseteq x^{n}D_{P}$. Therefore for each $P\in\starmax(D)$, $D_{P}$ is an AV-domain.
\end{proof}
Next we point out that a $\ast$-almost super-SH domain of type 2 is precisely a $\ast$-Krull domain.
\begin{corollary}
 Let $D$ be a domain and $\ast$ a finite character star-operation on $D$. The following statements are equivalent for $D$.
 \begin{itemize}
\item[(1)] $D$ is a $\ast$-almost super-SH domain of type 2.
\item[(2)] $D$ is a $\ast$-Krull domain.
\item[(3)] If $A$ is a finitely generated $\ast$-invertible ideal of $D$ with $A_{\ast}\neq D$, then $A_{\ast}$ is a $\ast$-product of $\ast$-almost super-homogeneous ideals of type 2.
\item[(4)] $D$ is a $\ast$-SH domain of type 2.
\end{itemize}
\end{corollary}
\begin{proof}
$(1)\Rightarrow (4)$\ \ Trivial.

$(4)\Leftrightarrow (2)$ \ \ \cite[Theorem 8]{AZ}.

$(2)\Rightarrow (3)$\ \ Suppose that $A$ is a finitely generated $\ast$-invertible $\ast$-homogeneous ideal of $D$ with $A_{\ast}\neq D$. Then $A_{\ast}=(A_{1} \cdots A_{k})_{\ast}$  by \cite[Theorem 8]{AZ}, where each $A_{i}$ is  a $\ast$-invertible $\ast$-homogeneous ideal of type 2 . And since a $\ast$-Krull domain is a $\ast$-almost IRKT, each $A_{i}$ is $\ast$-almost super-homogeneous by Theorem \ref{011}. Hence $A_{\ast}$ is a $\ast$-product of $\ast$-almost super-homogeneous ideals of type 2.

$(3)\Rightarrow (1)$\ \ Clear.
\end{proof}
\section{$\ast$-almost generalized Krull domains}
In this section we introduce $\ast$-AGKDs and we prove that a domain $D$ is a $\ast$-AGKD if and only if $D$ is a type 1 $\ast$-almost super-SH domain.
\begin{definition}\label{012}
\rm {Let $\ast$ be a finite character star-operation on the domain $D$. Then $D$ is called a \emph{$\ast$-almost generalized Krull domain} (\emph{$\ast$-AGKD}) if $D$ satisfies the following three conditions:
\begin{itemize}
\item[(1)] $D=\bigcap\{D_{P}\mid P\in X^{(1)}(D)\}$ is locally finite,
\item[(2)] $\starmax(D)=X^{(1)}(D)$, and
\item[(3)] $D_{P}$ is an AV-domain for each $P\in \starmax(D)$.
\end{itemize}}
\end{definition}

Let $D$ be a domain. If $D=\bigcap\{D_{P}\mid P\in X^{(1)}(D)\}$ is locally finite and for each $P\in X^{(1)}(D)$, $D_{P}$ is an AV-domain, then $D$ is said to be an \textit{AGKD}. Hence a domain $D$ is a $\ast$-AGKD if and only if $D$ is an AGKD and $\starmax(D)=X^{(1)}(D)$. Thus a $t$-AGKD is the same thing as an AGKD (see Proposition \ref{013}), while a $d$-AGKD is just a one-dimension finite character AP-domain.

Let $\Ass(K/D)$ be the set of associated primes of principal ideals of a domain $D$, $\Ass(K/D)=\{P\in \Spec(D)\mid P$ is minimal over $(aD:bD)$ for some $a,b\in D\}$. Then $\Ass(K/D)$ is a defining family for $D$ by \cite[Theorem E(i)]{T}. The function $g$ defined for all $A\in F(D)$ by $A\rightarrow\bigcap\{A_{P}\mid P\in \Ass(K/D)\}$ is a star-operation, which is so-called the \emph{$g$-operation}. This star-operation is also called the $\rho$-operation in \cite{Z01}. In particular, $g_{s}$ is a star-operation called the \emph{$f$-operation} in \cite{X}. Recall from \cite{QW02}, a domain $D$ is called a \emph{GW-domain} if $A_{g}=A_{w}$ for any $0\neq A\in F(D)$.  It is shown in \cite[Theorem 1.5]{QW02} that $D$ is a GW-domain if and only if ($\tmax(D)=$) $\wmax(D)\subseteq \Ass(K/D)$, if and only if $g$ is of finite character, i.e., $g=f$. Next, we point out that an AGKD is a GW-domain.
\begin{proposition}\label{013}
Let $D$ be a domain and $\ast$ a finite character star-operation on $D$.
\begin{itemize}
\item[(1)]If $X^{(1)}(D)=\starmax(D)$,  then $X^{(1)}(D)=\Ass(K/D)=\tmax(D)=\wmax(D)$.
\item[(2)]If $D$ is a $\ast$-AGKD, then $g=\ast_{w}=f=w$.
\end{itemize}
\end{proposition}

\begin{proof}\ \ $(1)$\ \ We first prove that $\Ass(K/D)= X^{(1)}(D)$. It is clear that $X^{(1)}(D)\subseteq \Ass(K/D)$. For the reverse inclusion, let $P\in \Ass(K/D)$. Then $P$ is a minimal prime ideal over $\ann_{D}(\frac{b}{a}+D)$ for some $\frac{b}{a}\in K$, where $K$ is the quotient field of $D$. Since $\ann_{D}(\frac{b}{a}+D)=(aD:_{R}bD)$ is a $v$-ideal of $D$, $P$ is a $t$-ideal of $D$  by \cite[Proposition 1.1 (5)]{HH}. Hence $P$ is a $\ast_{w}$-ideal of $D$.
Also, since $X^{(1)}(D)=\starmax(D)$, it follows from \cite[Corollary 2.10]{AC} that $\ast_{w}$ is precisely induced by $\{D_{Q}\mid Q\in X^{1}(D)\}$.  Hence by \cite[Theorem 1 (5)]{A}, there exists $Q\in X^{(1)}(D)$ such that $P\subseteq Q$.
But as $\height Q=1$, this forces $P=Q\in X^{(1)}(D)$. Hence $\Ass(K/D)\subseteq X^{(1)}(D)$, and so $X^{(1)}(D)=\starmax(D)=\Ass(K/D)$.
Now by \cite[Theorem E(i)]{T}, we have $D=\bigcap\{D_{P}\mid P\in \Ass(K/D)=\starmax(D)\}$. Thus $g=\ast_{w}$ has finite character. On other hand, since $X^{(1)}(D)=\Ass(K/D)$,  we have $\Ass(K/D)=\tmax{D}=\wmax{D}$ by \cite[Lemma 5.1 (5)(d)]{HZ}. Consequently, $X^{(1)}(D)=\Ass(K/D)=\tmax(D)=\wmax(D)$.

$(2)$\ \ By definition we have $X^{(1)}(D)=\starmax(D)$. Hence $X^{(1)}(D)=\Ass(K/D)=\tmax(D)=\wmax(D)$ by (1). So $g=\ast_{w}=w=f$ by \cite[Theorem 5]{QW02}.
\end{proof}
Let $\ast$ be a finite character star-operation on a domain $D$. Recall from \cite{HZ} that $P\in \starmax(D)$ is called \emph{$\ast$-potent} if it contains a $\ast$-homogeneous ideal and $D$ is called \emph{$\ast$-potent} if each $P\in \starmax(D)$ is $\ast$-potent.
\begin{theorem}\label{014}
Let $D$ be a domain and $\ast$ a finite character star-operation on $D$. Then the following statements are equivalent for $D$.
\begin{itemize}
\item[(1)] $D$ is a $\ast$-AGKD.
\item[(2)] $D$ is an AP$\ast$MD and a $\ast$-WKD.
\item[(3)] $D$ is a $\ast$-potent AP$\ast$MD with $X^{(1)}(D)=\starmax(D)$.
\item[(4)] $D$ is an $\ast$-almost IRKT and a $\ast$-WKD.
\item[(5)] $D$ is an $\ast$-almost IRKT and every $\ast$-invertible $\ast$-homogeneous ideal has type 1.
\item[(6)] $D$ is a $\ast$-h-local domain and every $\ast$-invertible $\ast$-homogeneous ideal is $\ast$-almost super-homogeneous and has type 1.
\item[(7)] $D$ is a type 1 $\ast$-almost super-SH domain.
\item[(8)] If $A$ is a finitely generated $\ast$-invertible ideal of $D$ with $A_{\ast}\neq D$, then $A_{\ast}$ is a $\ast$-product of type 1 $\ast$-almost super-homogeneous ideals.

\end{itemize}
\end{theorem}
\begin{proof}

$(1)\Rightarrow (2)$ \ \ This follows from the definitions and \cite[Theorem 2.4]{L}.

$(2)\Rightarrow (3)$\ \ Since $D$ is a $\ast$-WKD, $D$ is a $\ast$-h-local domain and $X^{(1)}(D)=\starmax(D)$. Hence $D$ is a $\ast$-SH domain by \cite[Theorem 4]{AZ}. Let $P\in \starmax(D)$. Take $0\neq x\in P$. Then $(A_{1}\cdots A_{k})_{\ast}=xD$ where each $A_{i}$ is a $\ast$-homogeneous ideal of $D$. Hence $A_{1}\cdots A_{k}\subseteq P$. It follows that there exists some $A_{j}$ such that $A_{j}\subseteq P$. Hence $P$ is $\ast$-potent. So $D$ is a $\ast$-potent domain.

$(3)\Rightarrow (4)$\ \ Since $X^{(1)}(D)=\starmax(D)$ by (3),  we have $\Ass(K/D)=X^{(1)}(D)=\starmax(D)$ by Proposition \ref{013}. Hence $D=\mathop{\bigcap}\{D_{P}\mid P\in X^{(1)}(D)\}$  by \cite[Theorem E(ii)]{T}.  Also, since $D$ is a $\ast$-potent domain with $X^{(1)}(D)=\starmax(D)$, $D$ has finite $\ast$-character by \cite[Theorem 5.3]{HZ}. Hence $D=\mathop{\bigcap}\{D_{P}\mid P\in X^{(1)}(D)\}$ is locally finite. It follows that $D$ is a WKD.  So $D$ is a $\ast$-h-local domain by \cite[Theorem 7]{AZ}. On other hand, since $D$ is an AP$\ast$MD, $D_{P}$ is an AV-domain for each $P\in \starmax(D)$ by \cite[Theorem 2.4]{L}. So $D$ is a $\ast$-almost IRKT.

$(4)\Rightarrow (1)$ \ \ Obvious.

$(4)\Rightarrow (5)$ \ \ Since $D$ is a $\ast$-WKD, every $\ast$-homogenous ideal of $D$ has type 1 by \cite[Theorem 7]{AZ}. Hence every $\ast$-invertible $\ast$-homogeneous ideal of $D$ has type 1.

$(5)\Rightarrow (6)$\ \ This follows from Theorem \ref{011}.

$(6)\Rightarrow(7)$\ \ Since $D$ is a $\ast$-h-local domain, it follows from \cite[Theorem 4]{AZ} that $D$ is a $\ast$-SH domain. Let $xD$ be a nonzero proper principal ideal of $D$. Then $xD=(A_{1}\cdots A_{k})_{\ast}$ where each $A_{i}$ is a $\ast$-invertible $\ast$-homogeneous ideal of $D$. Hence $A_{i}$ is  $\ast$-almost super-homogeneous and has type 1 by (4). So $D$ is a type 1 $\ast$-almost super-SH domain.

$(7)\Rightarrow (4)$\ \ This follows from Theorem \ref{011} and \cite[Theorem 7]{AZ}.

$(6)\Rightarrow (8)$\ \ Since $D$ is $\ast$-h-local, $D$ is a $\ast$-SH domain. Hence by \cite[Theorem 6]{AZ} $A_{\ast}=(A_{1}\cdots A_{k})_{\ast}$ where each $A_{i}$ is $\ast$-homogeneous. Since $A$ is $\ast$-invertible, each $A_{i}$ is also $\ast$-invertible. Hence $A_{i}$ is $\ast$-almost super-homogeneous and has type 1  by (6). It follows that $A$ is a $\ast$-product of type 1 $\ast$-almost super-homogeneous ideals.

$(8)\Rightarrow (7)$\ \ This is clear.
\end{proof}

Now we give an example to show that (1) a $\ast$-AGKD is not necessarily a $\ast$-GKD, (2) a $\ast$-almost IRKT is not necessarily a $\ast$-IRKT, and (3) a $\ast$-almost super-homogeneous ideal is not necessarily a $\ast$-super-homogeneous ideal.
\begin{example}\label{21}
\rm{ Let $\mathbb{F}$ be a field of characteristic 2. Set $D=\mathbb{F}[X^{2},X^{3}]=\mathbb{F}+X^{2}\mathbb{F}[X]$. Then $D$ is an AGCD-domain by \cite[Corollary 3.2]{ACP1}, and $D$ is a WKD by \cite[Corollary 4.6]{ACP}. Hence $D$ is a WKD and an AP$v$MD by \cite[Theorem 3.1]{L02}. So $D$ is an AGKD ($t$-AGKD) by Theorem \ref{014}. On the other hand, since $D$ is never integrally closed, $D$ is not a P$v$MD. Hence $D$ is not a GKD by \cite[Theorem 5.9]{HZ}. So $D$ is an AGKD, but not a GKD. By Theorem \ref{014} and \cite[Theorem 11]{AZ} it follows that $D$ is an almost IRKT, but not a IRKT. Also, since $D$ is an AGKD but not a GKD, $D$ is an $t$-almost super-SH domain but not a $t$-super-SH domain. Hence there exists a $t$-almost super-homogenous ideal of $D$ but not a $t$-super-homogeneous ideal.}

\end{example}

\section{$\ast$-almost factorial general-SH domains}
In this section we introduce $\ast$-almost factorial general-SH domains. We prove that a domain $D$ is a $\ast$-almost factorial general-SH domain if and only if $D$ is a $\ast$-IRKT and an AGCD.
\begin{definition}\label{015}
\rm{Let $\ast$ be a finite character star-operation on a domain $D$.
\begin{itemize}
\item[(1)]A ideal $A$ of $D$ is called \emph{$\ast$-almost factorial general-homogeneous} (\emph{$\ast$-afg-homogeneous}) if
\begin{itemize}
\item[(i)] $A$ is a $\ast$-invertible $\ast$-homogeneous ideal, and
\item[(ii)] given $b_{1},\dots, b_{s}\in P$ with $A^{r}\subseteq (b_{1}, \dots, b_{s})_{\ast}$ for some $r\in \mathbb{N}$, there exists an $n\in \mathbb{N}$ with $(b_{1}^{n},\dots, b_{s}^{n})$ principal.
\end{itemize}
\item[(2)] A domain $D$ is called a \emph{$\ast$-almost factorial general-SH domain} (\emph{$\ast$-afg-SH domain}) if every nonzero proper principal ideal of $D$ is a $\ast$-product of $\ast$-afg-homogeneous ideals.
\end{itemize}}
\end{definition}
\remark\label{24}
In Definition \ref {015} (1), the $n$ depends on $b_{1},\dots, b_{s}$.

\begin{proposition}\label{016}
Let $\ast$ be a finite character star-operation on a domain $D$ and $A$ a $P$-$\ast$-afg-homogeneous ideal of $D$. Then the following statements hold.
\begin{itemize}
\item[(1)] $A$ is  $\ast$-almost super-homogeneous.

\item[(2)] $(A^{n})_{\ast}$ is a principal ideal of $D$ for some $n\in \mathbb{N}$.

\item[(3)] If $(b_{1},\dots, b_{s})$ is a $P$-$\ast$-homogeneous ideal of $D$, then $A^{n}\subseteq (b_{1}^{n},\dots,b_{s}^{n})_{\ast}$ or $(b_{1}^{n},\dots,b_{s}^{n})\subseteq (A^{n})_{\ast}$ for some $n\in \mathbb{N}$.

\item[(4)] If $B$ is a $P$-$\ast$-afg-homogeneous ideal of $D$, then $A^{n}\subseteq (B^{n})_{\ast}$ or $B^{n}\subseteq (A^{n})_{\ast}$ for some $n\in \mathbb{N}$.

\item[(5)] If $B$ is a $P$-$\ast$-afg-homogeneous ideal of $D$, so is $AB$.

\item[(6)] $A^{n}$ is  $P$-$\ast$-afg-homogeneous for each $n \in \mathbb{N}$.
\end{itemize}
\end{proposition}
\begin{proof}
 $(1)$\ \ This follows from the definitions.

$(2)$\ \  Let $A=(a_{1},\dots, a_{k})$. Then $(a_{1}^{n},\dots, a_{k}^{n})_{\ast}$ is a principal ideal of $D$ for some $n\in \mathbb{N}$. Since $A$ is $\ast$-invertible, we have $(A^{n})_{\ast}=(a_{1}^{n},\dots, a_{k}^{n})_{\ast}$ by \cite[Lemma 2.2]{L}. Hence $(A^{n})_{\ast}$ is a principal ideal of $D$.

$(3)$\ \ Since $A$ is $\ast$-afg-homogeneous, $A$ is $\ast$-almost super-homogeneous by (1). Hence $A^{n}\subseteq (b_{1}^{n},\dots, b_{s}^{n})_{\ast}$ or $(b_{1}^{n},\dots, b_{s}^{n})\subseteq (A^{n})_{\ast}$ by Proposition \ref{005}(1).

$(4)$\ \ Suppose that $B=(b_{1},\dots,b_{s})$. Then by (3) $A^{n}\subseteq (b_{1}^{n},\dots, b_{s}^{n})_{\ast}$ or $(b_{1}^{n},\dots, b_{s}^{n})$ $\subseteq (A^{n})_{\ast}$ for some $n\in \mathbb{N}$. Since $B$ is $\ast$-afg-homogeneous, $B$ is $\ast$-invertible. Hence, $(B^{n})_{\ast}=(b_{1}^{n},\dots, b_{s}^{n})_{\ast}$ by \cite[Lemma 2.2]{L}. So $A^{n}\subseteq (B^{n})_{\ast}$ or $B^{n}\subseteq (A^{n})_{\ast}$.

$(5)$\ \ By (4) we have $A^{n}\subseteq (B^{n})_{\ast}$ or $B^{n}\subseteq (A^{n})_{\ast}$ for some $n\in \mathbb{N}$. Let $C=(c_{1},\dots, c_{l})$ be $P$-$\ast$-homogeneous with $(AB)^{r}\subseteq (c_{1},\dots, c_{l})_{\ast}$ for some $r\in \mathbb{N}$. Then by Lemma \ref{002} $C^{nl}\subseteq (c_{1}^{n},\dots, c_{l}^{n})$. Hence $A^{nrl}B^{nrl}\subseteq (A^{nrl}B^{nrl})_{\ast}\subseteq (C^{nl})_{\ast}\subseteq (c_{1}^{n},\dots, c_{l}^{n})_{\ast}$. If $A^{n}\subseteq (B^{n})_{\ast}$, then $A^{2nrl}\subseteq (c_{1}^{n},\dots, c_{l}^{n})_{\ast}$. Since $A$ is $\ast$-afg-homogenous, $(c_{1}^{mn},\dots, c_{l}^{mn})_{\ast}$ is principal for some $m\in \mathbb{N}$. Similarly if $B^{n}\subseteq (A^{n})_{\ast}$, then $(c_{1}^{sn},\dots, c_{l}^{sn})_{\ast}$ is principal for some $s\in \mathbb{N}$. Also, it is clear that $AB$ is $\ast$-invertible and similar to both $A$ and $B$ by \cite[Proposition 2]{AZ}. Hence $AB$ is $P$-$\ast$-afg-homogeneous.

$(6)$\ \ This follows from $(5)$.

\end{proof}

\begin{corollary}\label{017}
If $D$ is a $\ast$-afg-SH domain, then $D$ is a $\ast$-almost super-SH domain.
\end{corollary}
\begin{proof}
This follows from Proposition \ref{016}.
\end{proof}
Now by Proposition \ref{016} (5) a product of similar $\ast$-afg-homogeneous
ideals is again $\ast$-afg-homogeneous. Thus the proof of Theorem \ref{009}
gives the corresponding uniqueness result for $\ast$-products of $\ast$-afg-homogeneous
ideals.
\begin{theorem}\label{0009}
Let $\ast$ be a finite character star-operation on a domain $D$ and let $A_{1},\dots, A_{n}$ be $\ast$-afg-homogeneous ideals of $D$. Then the $\ast$-product $(A_{1}\cdots A_{n})_{\ast}$ can be expressed uniquely, up to order, as a product of pairwise $\ast$-comaximal $\ast$-afg-homogeneous ideals.
\end{theorem}

Let $\ast$ be a finite character star-operation on a domain $D$. The set $\ast\mbox{-Inv}(D)$ of $\ast$-invertible
fractional $\ast$-ideals forms a group under the $\ast$-product $I\ast J :=(IJ)_{\ast}$ with subgroup
$\mbox{Prin}(D)$, the set of nonzero principal fractional ideals of $D$. The quotient group
$Cl_{\ast}(D):=\ast\mbox{-Inv}(D)/\mbox{Prin}(D)$ is called the \emph{$\ast$-class group} of $D$ in \cite{AF}. If $\ast_{1}\leq \ast_{2}$ are finite character star-operations on $D$, then $Cl_{\ast_{1}}(D)\subseteq Cl_{\ast_{2}}(D)$. Let $\ast$ be a finite character star-operation on a domain $D$.
Then $D$ is called a \emph{$\ast$-almost B\'{e}zout domain} in \cite{AZ} if for $0\neq a, b \in D$, there exists an $n=n(a,b)\in \mathbb{N}$ with $(a^{n},b^{n})_{\ast}$ principal. It follows from \cite[Theorem 3.4]{L} that a domain $D$ is a $\ast$-almost B\'{e}zout domain if and only if $D$ is an AP$\ast$MD with $Cl_{\ast}(D)$ torsion. If $\ast_{1}\leq \ast_{2}$ are finite character star-operations on $D$, then $D$ $\ast_{1}$-almost B\'{e}zout implies $D$ is a $\ast_{2}$-almost B\'{e}zout. Next we characterize $\ast$-afg-SH domains and we need the following lemma.
\begin{lemma}\label{018}
Let $\ast$ be a finite character star-operation on a domain $D$. If $D$ is a $\ast$-almost IRKT, then $D$ is an AGCD-domain if and only if $D$ is a $\ast$-almost B\'{e}zout domain.
\end{lemma}

\begin{proof}
($\Rightarrow$) Let $D$ be an AGCD-domain. Then $Cl_{t}(D)$ is torsion by \cite[Theorem 3.1]{L02}. Since $Cl_{\ast}(D)\subseteq Cl_{t}(D)$, it follows that $Cl_{\ast}(D)$ is torsion. Also since $D$ is a $\ast$-almost IRKT, $D$ is a AP$\ast$MD. Hence $D$ is a $\ast$-almost B\'{e}zout domain by \cite[Theorem 3.4]{L}.

($\Leftarrow$) Let $D$ be a $\ast$-almost B\'{e}zout domain. Then $D$ is a $t$-almost B\'{e}zout domain and hence $D$ is an AGCD-domain.
\end{proof}

\begin{theorem}\label{019}
 Let $D$ be a domain and $\ast$ a finite character star-operation on $D$. The following statements are equivalent for $D$.
\begin{itemize}
\item[(1)] $D$ is a $\ast$-afg-SH domain.
\item[(2)] $D$ is a $\ast$-almost IRKT and an AGCD-domain.
\item[(3)] $D$ is a $\ast$-almost IRKT with $Cl_{\ast}(D)$ torsion.
\item[(4)] $D$ is a $\ast$-h-local domain and each $\ast$-invertible $\ast$-homogeneous ideal is $\ast$-afg-homogeneous.
\end{itemize}
\end{theorem}
\begin{proof} $(1)\Rightarrow (2)$\ \ Since $D$ is a $\ast$-afg-SH domain, $D$ is a $\ast$-almost super-SH domain by Corollary \ref{017}.  Hence $D$ is a $\ast$-almost IRKT by Theorem \ref{011}. We only need to show that $D$ is an AGCD-domain. Let $c$ be nonzero nonunit of $D$. Then $cD=(C_{1}\cdots C_{l})_{\ast}$, where the $C_{i}$ are mutually
$\ast$-comaximal and $\ast$-afg homogeneous. Hence by Proposition \ref{016} (2), $(C_{i}^{n_{i}})_{\ast}=c_{i}'D$ for some $n_{i}\in \mathbb {N}$ and $c_{i}'\in D$ ($i=1,\dots, l$). Set $n=\prod_{i=1}^{l}n_{i}$ and $c_{i}=(c_{i}')^{n/n_{i}}$. Then $(C_{i}^{n})_{\ast}=c_{i}D$ is $M(C_{i})$-$\ast$-homogeneous. Hence $c^{n}D=(c_{1}D\cdots c_{l}D)_{\ast}=c_{1}D\cdots c_{l}D$, and the $c_{i}D$ are mutually
$\ast$-comaximal. So $c^{n}D=c_{1}D\cdots c_{l}D=c_{1}D\bigcap \cdots\bigcap c_{l}D$. Let $a$ and $b$ be nonzero nonunits of $D$. Let $P_{1},\dots, P_{s}$ be the maximal $\ast$-ideal containing $a$ or $b$. Then for suitable $m$ $a^{m}D=a_{1}D\bigcap\cdots \bigcap a_{s}D=a_{1}D\cdots a_{s}D$ and $b^{m}D=b_{1}D\bigcap\cdots \bigcap b_{s}D=a_{1}D\cdots b_{s}D$ where either $a_{i}D$ (resp., $b_{i}D$) is $P_{i}$-$\ast$-homogeneous or $a_{i}=1$ (resp., $b_{i}=1$).
Thus $a_{i}^{m_{i}}D\subseteq b_{i}^{m_{i}}D$ or $a_{i}^{m_{i}}D\supseteq b_{i}^{m_{i}}D$ for some $m_{i}\in \mathbb{N}$ ($i=1,\dots, s$) by Proposition \ref{016}(4). It follows that $a_{i}^{m_{i}}D\bigcap b_{i}^{m_{i}}D$ is principal and $P_{i}$-homogeneous ($i=1,\dots, s$). Set $n=\prod_{i} m_{i}$ and $n_{i}=\prod_{j\neq i} m_{j}$. Then $a^{mn}D=(a_{1}^{m_{1}}D)^{n_{1}}\bigcap \cdots\bigcap (a_{s}^{m_{s}}D)^{n_{s}}$ and  $b^{mn}D=(b_{1}^{m_{1}}D)^{n_{1}}\bigcap \cdots\bigcap (b_{s}^{m_{s}}D)^{n_{s}}$. Hence $a^{mn}D\bigcap b^{mn}D=(a_{1}^{m_{1}n_{1}}D\bigcap b_{1}^{m_{1}n_{1}}D)\bigcap \cdots \bigcap $ $(a _{s}^{m_{s}n_{s}}D)\bigcap (b_{s}^{m_{s}n_{s}}D)$ and $(a_{i}^{m_{i}}D)^{n_{i}}\bigcap (b_{i}^{m_{i}}D)^{n_{i}}$ is principal and $P_{i}$-$\ast$-homogeneous ($i=1,\dots, r$). So $a^{mn}D\bigcap b^{mn}D$ is a product of principal ideals and hence is a principal ideal. Consequently, $D$ is an AGCD-domain.

$(2)\Rightarrow (3)$\ \ Since $D$ is an AGCD-domain and a $\ast$-almost IRKT by (2), $D$ is a $\ast$-almost B\'{e}zout domain by Lemma \ref{018}. Hence $Cl_{\ast}(D)$ is torsion by \cite[Theorem 3.1]{L02}.

$(3)\Rightarrow (4)$\ \ Since $D$ is a $\ast$-almost IRKT, it is clear that $D$ is a $\ast$-h-local domain. Now suppose that $A$ is a $\ast$-invertible and $P$-$\ast$-homogeneous ideal of $D$. Then $A$ is $\ast$-almost super homogeneous by Theorem \ref{011}. Let $(b_{1},\dots, b_{s})$ be $P$-$\ast$-homogeneous with $A^{r}\subseteq (b_{1},\dots, b_{s})_{\ast}$ for some $r\in\mathbb{N}$. Then $(b_{1}^{n},\dots, b_{s}^{n})$ is $\ast$-invertible for some $n\in \mathbb{N}$.  Hence $((b_{1}^{n},\cdots, b_{s}^{n})^{m})_{\ast}$ is principal for some $m\in \mathbb{N}$ because $Cl_{\ast}(D)$ is torsion. Since $(b_{1}^{n},\dots, b_{s}^{n})$ is $\ast$-invertible, $((b_{1}^{n},\dots, b_{s}^{n})^{m})_{\ast}=(b_{1}^{mn},\dots, b_{s}^{mn})_{\ast}$ by Lemma \ref{003}. It follows that $(b_{1}^{mn},\dots, b_{s}^{mn})_{\ast}$ is principal. So $A$ is $\ast$-afg-homogeneous.

$(4)\Rightarrow (1)$\ \ Clear.
\end{proof}
Recall from \cite{AZ} that a $\ast$-homogeneous ideal $A$ of a domain $D$ is called \emph{$\ast$-almost factorial-homogeneous ($\ast$-af-homogeneous)} if for each $\ast$-homogeneous ideal $B\supseteq A $, there exists some $n\in \mathbb{N}$ with $(B^{n})_{\ast}$ principal. The domain $D$ is a
\emph{$\ast$-af-SH domain} if for each nonzero nonunit $x\in D$, $xD$ is expressible as a $\ast$-product of finitely many $\ast$-af-homogeneous ideals. Next we point out that every $\ast$-af SH domain is a $\ast$-af-SH domain.
\begin{corollary}\label{020}
Every $\ast$-af-SH domain is a $\ast$-afg-SH domain.
\end{corollary}
\begin{proof}
Let $D$ be a $\ast$-fg-SH domain. Then $D$ is a $\ast$-IRKT and an AGCD-domain by \cite[Theorem 13]{AZ}. Hence $D$ is a $\ast$-almost IRKT and an AGCD-domain.  It follows from Theorem \ref{019} that $D$ is a $\ast$-afg-SH domain.
\end{proof}
\begin{corollary}\label{022}
 Let $D$ be a domain and $\ast$ a finite character star-operation on $D$. The following statements are equivalent for $D$.
 \begin{itemize}
\item[(1)] $D$ is a $\ast$-afg-SH domain of type 1.
\item[(2)] $D$ is a $\ast$-AGKD and an AGCD-domain.
\item[(3)] $D$ is a $\ast$-h-local domain and each $\ast$-invertible $\ast$-homogeneous ideal is a $\ast$-afg-homogeneous ideal of type 1.
\item[(4)] If $A$ is a finitely generated $\ast$-invertible ideal of $D$ with $A_{\ast}\neq D$, then $A_{\ast}$ is a $\ast$-product of $\ast$-afg-homogeneous ideals of type 1.
\item[(5)] $D$ is a $\ast$-AGKD with $Cl_{\ast}(D)$ torsion.
\end{itemize}
\end{corollary}
\begin{proof}
$(1)\Rightarrow (2)$\ \ Since $D$ is a $\ast$-afg-SH domain of type 1, $D$ is a $\ast$-almost super-SH domain of type 1 by Corollary \ref{017}. Hence $D$ is a $\ast$-AGKD by Theorem \ref{014}. Also, since $D$ is a $\ast$-afg-SH domain, $D$ is an AGCD-domain by Theorem \ref{019}. So $D$ is a $\ast$-AGKD and an AGCD-domain.

$(2)\Rightarrow (3)$\ \ Since $D$ is a $\ast$-AGKD and an AGCD-domain, $D$ is a $\ast$-almost IRKT and an AGCD-domain. Hence by Theorem \ref{019} $D$ is a $\ast$-h-local domain and each $\ast$-invertible $\ast$-homogeneous ideal is $\ast$-afg-homogeneous. By Theorem \ref{014} it follows that each $\ast$-invertible $\ast$-homogeneous ideal has type 1. So $D$ is a $\ast$-SH domain and each $\ast$-invertible $\ast$-homogeneous ideal is a $\ast$-afg-homogeneous ideal of type 1.

$(3)\Rightarrow (4)$\ \ Suppose that $A$ is a finitely generated $\ast$-invertible ideal of $D$ with $A_{\ast}\neq D$. Then by \cite[Theorem 6]{AZ} $A=(A_{1}\cdots A_{k})$ where each $A_{i}$ is $\ast$-homogeneous. Since $A$ is $\ast$-invertible, each $A_{i}$ is $\ast$-invertible. Hence each $A_{i}$ is $\ast$-afg-homogeneous ideals of type 1.

$(4)\Rightarrow (1)$\ \ Clear.

$(2)\Leftrightarrow (5)$\ \ This follows from Theorem \ref{019}.
\end{proof}
\begin{corollary}\label{023}
 Let $D$ be a domain and $\ast$ a finite character star-operation on $D$. The following statements are equivalent for $D$.
 \begin{itemize}
\item[(1)] $D$ is a $\ast$-afg-SH domain of type 2.
\item[(2)] $D$ is a $\ast$-Krull domain and an AGCD-domain.
\item[(3)] If $A$ is a finitely generated $\ast$-invertible ideal of $D$ with $A_{\ast}\neq D$, then $A_{\ast}$ is a $\ast$-product of $\ast$-afg-homogeneous ideals of type 2.
\item[(4)] $D$ is a $\ast$-Krull domain with $Cl_{\ast}(D)$ torsion.
\item[(5)] $D$ is a $\ast$-af-SH domain.
\end{itemize}
\end{corollary}
\begin{proof}$(1)\Rightarrow (5)$\ \ Trivial.

$(5)\Leftrightarrow (4)\Leftrightarrow (2)$\ \ \cite[Theorem 15]{AZ}.

$(2)\Rightarrow (3)$ \ \ Suppose that $A$ is a finitely generated $\ast$-invertible ideal of $D$ with $A_{\ast}\neq D$. Then  $A_{\ast}=(A_{1} \cdots A_{k})_{\ast}$  by \cite[Theorem 8]{AZ}, where each $A_{i}$ is a $\ast$-invertible $\ast$-homogeneous ideal of type 2 . Also since $D$ is a $\ast$-almost IRKT and an AGCD-domain, each $A_{i}$ is $\ast$-afg-homogeneous. Hence $A_{\ast}$ is a $\ast$-product of $\ast$-afg-homogeneous ideals of type 2.

$(3)\Rightarrow (1)$ \ \ Clear.
\end{proof}

Finally, we give an example to show that (1) a $\ast$-afg-SH domain is not necessarily a $\ast$-fg-SH domain and (2) a $\ast$-afg homogeneous ideal is not necessarily a $\ast$-fg homogeneous ideal.
\begin{example}\label{021}
\rm{ Let $K$ be a field of characteristic $p > 0$ and let $K\subset L$ be
a purely inseparable field extension. Set $D=K+XL[X]$. Then $D$ is a WKD by \cite[Corollary 3.11]{ACP}. Hence $D$ is a $t$-h-local domain. Since $D$ is an AB-domain by \cite[Example 4.14]{AZ02}, $D$ is a AGCD-domain. Hence $D$ is a AP$v$MD by \cite[Theorem 3.1]{L02}. So $D_{P}$ is an AV-domain for each $P\in\tmax(D)$ by \cite[Theorem 2.3]{L02}. It follows that $D$ is a $t$-almost IRKT. Thus $D$ is a $t$-afg-SH domain by Theorem \ref{019}. Also since the integral closure of $D$ is precisely $L[X]$, it follows that $D$ is not integrally closed. Hence $D$ is not a $t$-IRKT. So $D$ is not a $t$-fg-SH domain by \cite[Theorem 13]{AZ}. Thus, there exists a $t$-afg homogeneous ideal in $D$, that is a $t$-fg homogeneous ideal.}

\end{example}

\end{document}